**Hankel determinants of generalized q-exponential polynomials**

Johann Cigler

Fakultät für Mathematik, Universität Wien

johann.cigler@univie.ac.at

**Abstract**

Recently I. Mezö studied a simple but interesting generalization of the exponential polynomials. In this note I consider two $q-$analogues of these polynomials and compute their Hankel determinants.

## 1. Introduction

In [5] I. Mezö studied a simple but interesting generalization of the exponential polynomials. In this note I consider two $q-$analogues of these polynomials and compute their Hankel determinants.

Let $r$ be a real number. Define (generalized) $q-$Stirling numbers (of the second kind) $S(n,k,r)$ by the recurrence

$$S(n,k,r)=S(n-1,k-1,r)+[k+r]S(n-1,k,r) \tag{1.1}$$

with boundary values $S(0,k,r)=\left[k=0\right]$ and $S(n,0,r)=[r]^n$.

As usual in $q-$analysis we write $[r]=\frac{1-q^r}{1-q}$ for $r\in\mathbb{R}$, $[n]!=\prod_{i=0}^{n}[i]$ for $n\in\mathbb{N}$, $\begin{bmatrix} n \\ k \end{bmatrix}=\frac{[n]!}{[k]![n-k]!}$ for $k\in\mathbb{N}$ with $0\le k\le n$ and $\begin{bmatrix} n \\ k \end{bmatrix}=0$ else. By $D$ we denote the $q-$differentiation operator defined by $Df(x)=\frac{f(x)-f(qx)}{x-qx}$, which satisfies $Dx^r=[r]x^{r-1}$.

The generating function of these $q-$Stirling numbers is given by

$$\sum_{n\ge 0}S(n,k,r)z^n=\frac{z^k}{(1-[r]z)(1-[r+1]z)\cdots(1-[r+k]z)}. \tag{1.2}$$

Let

$$\langle x\rangle_{r,k}=\prod_{j=0}^{k-1}\left(x-[r+j]\right). \tag{1.3}$$

Then it is easily verified that

$$\sum_{k=0}^{n} S(n,k,r)\langle x\rangle_{r,k} = x^n. \tag{1.4}$$

We introduce now two kinds of generalized $q-$exponential polynomials

$$\varphi_n(x,r) = \sum_{k=0}^{n} S(n,k,r)x^k \tag{1.5}$$

and

$$\Phi_n(x,r) = \sum_{k=0}^{n} S(n,k,r)q^{\binom{k}{2}}\left(q^r x\right)^k. \tag{1.6}$$

For $r=0$ these are the usual $q-$analogues of the exponential polynomials.

I want to compute their Hankel determinants.

I shall use the following facts (cf. e.g. [2],[3] or [4]) :

Let $(a_n)$ be a sequence such that $a_0 = 1.$

Let $d(n,k) = \det\left(a_{i+j+k}\right)_{i,j=0}^{n-1}$ denote their Hankel determinants.

Define a linear functional $F$ on the polynomials by $F(x^n) = a_n.$

Suppose there exists a sequence of orthogonal polynomials $p_n(x) = x^n + c_1 x^{n-1} + \cdots + c_n$ with respect to $F.$ This means that $F(p_n p_k) = d_n[n=k]$ with $d_n \neq 0.$

Then

$$d(n,0) = \prod_{i=0}^{n-1} d_i \tag{1.7}$$

and

$$d(n,1) = d(n,0)(-1)^n p_n(0). \tag{1.8}$$

## 2. The polynomials $\varphi_n(x,r)$.

**Theorem 2.1**

*The Hankel determinants of* $\varphi_n(x,r)$ *are given by*

$$d(n,0)=q^{\binom{n}{3}}\left(q^r x\right)^{\binom{n}{2}}\prod_{k=0}^{n-1}[k]! \tag{2.1}$$

*and*

$$d(n,1)=q^{\binom{n}{3}}\left(q^r x\right)^{\binom{n}{2}}\prod_{k=0}^{n-1}[k]!\sum_{k=0}^{n}\begin{bmatrix} n \\ k \end{bmatrix}q^{\binom{k}{2}}x^k\prod_{j=0}^{n-k-1}[r+j]. \tag{2.2}$$

In order to prove this we consider the linear operator $U_r$ on the polynomials defined by

$$U_r\left\langle x\right\rangle_{r,n}=x^n. \tag{2.3}$$

Then

$$U_r x U_r^{-1}=x\left(1+x^{-r}Dx^r\right). \tag{2.4}$$

For $U_r x U_r^{-1}x^n=U_r x\left\langle x\right\rangle_{r,n}=U_r\left(\left\langle x\right\rangle_{r,n+1}+[r+n]\left\langle x\right\rangle_{r,n}\right)=x^{n+1}+[r+n]x^n=x\left(1+x^{-r}Dx^r\right)x^n.$

Let $F_r$ be the linear functional defined by

$$F_r\left(\left\langle x\right\rangle_{r,n}\right)=a^n. \tag{2.5}$$

Then the orthogonal polynomials with respect to $F_r$ are given by

$$h_n(x,a,r)=\sum_{k=0}^{n}(-a)^k q^{\binom{k}{2}}\begin{bmatrix} n \\ k \end{bmatrix}\left\langle x\right\rangle_{r,n-k}. \tag{2.6}$$

This is a variant of the $q-$Poisson-Charlier polynomials.

They satisfy the recurrence relation

$$h_{n+1}(x,a,r)=\left(x-[n+r]-q^n a\right)h_n(x,a,r)-q^{r+n-1}a[n]h_{n-1}(x,a,r). \tag{2.7}$$

To prove this we consider $p_n(x,a)=\prod_{k=0}^{n-1}\left(x-q^k a\right)=\sum_{k=0}^{n}(-a)^k q^{\binom{k}{2}}\begin{bmatrix} n \\ k \end{bmatrix}x^{n-k}$. Then $Dp_n(x,a)=[n]p_{n-1}(x,a)$ (see e.g. [1]).

We have

$$U_r h_n(x,a,r) = \sum_{k=0}^{n} (-a)^k q^{\binom{k}{2}} \begin{bmatrix} n \\ k \end{bmatrix} x^{n-k} = p_n(x,a). \tag{2.8}$$

(2.4) implies

$$\begin{aligned} U_r x h_n(x,a,r) &= U_r x U_r^{-1} p_n(x,a) = x\left(1 + x^{-r} D x^r\right) p_n(x,a) \\ &= p_{n+1}(x,a) + q^n a p_n(x,a) + [r] p_n(x,a) + q^r x[n] p_{n-1}(x,a) \\ &= p_{n+1}(x,a) + q^n a p_n(x,a) + [r] p_n(x,a) + q^r [n] p_n(x,a) + q^{r+n-1}[n] a p_{n-1}(x,a) \end{aligned}$$

By applying $U_r^{-1}$ we get

$$x h_n(x,a,r) = h_{n+1}(x,a,r) + \left([r] + q^r[n] + q^n a\right) h_n(x,a,r) + q^{r+n-1}[n] a h_{n-1}(x,a). \tag{2.9}$$

It is clear that $F_r\left(h_n(x,a,r)\right) = p_n(a,a) = 0$ for $n > 0$.

By (2.9) this implies $F_r\left(x^k h_n(x,a,r)\right) = 0$ for $k < n$ and

$$F_r\left(x^n h_n(x,a,r)\right) = q^{r+n-1}[n] a F_r\left(x^{n-1} h_{n-1}(x,a,r)\right) = \prod_{k=1}^{n} q^{r+k-1}[k] a = \left(q^r a\right)^n q^{\binom{n}{2}} [n]! \tag{2.10}$$

This implies (2.1).

For the second Hankel determinant we note that

$$(-1)^n h_n(0,a,r) = \sum_{k=0}^{n} (-a)^k q^{\binom{k}{2}} \begin{bmatrix} n \\ k \end{bmatrix} \langle 0 \rangle_{r,n-k} = \sum_{k=0}^{n} a^k q^{\binom{k}{2}} \begin{bmatrix} n \\ k \end{bmatrix} [r][r+1]\cdots[r+n-k-1].$$

**Remark**

For the usual $q-$Stirling numbers $S(n,k) = S(n,k,0)$ the formula

$$(q-1)^{n-k} S(n,k) = \sum_i (-1)^{n-i} \binom{n}{i} \begin{bmatrix} i \\ k \end{bmatrix} \tag{2.11}$$

holds. This can be generalized to

$$(q-1)^{n-k} q^{rk} S(n,k,r) = \sum_i (-1)^{n-i} q^{ri} \binom{n}{i} \begin{bmatrix} i \\ j \end{bmatrix}. \tag{2.12}$$

For by changing $x \to \frac{x-1}{q-1}$ it is easily verified that

$$\sum_{k=0}^{n} \begin{bmatrix} n \\ k \end{bmatrix} (q-1)^k q^{r(n-k)} \langle x \rangle_{r,k} = \left(1-(1-q)x\right)^n = \sum_{k=0}^{n} \binom{n}{k} (q-1)^k x^k. \qquad (2.13)$$

Applying $U_r \langle x \rangle_{r,n} = x^n$ we get

$$\sum_{k=0}^{n} \begin{bmatrix} n \\ k \end{bmatrix} (q-1)^k q^{r(n-k)} x^k = \sum_{k=0}^{n} \binom{n}{k} (q-1)^k \varphi_k(x,r). \qquad (2.14)$$

This is equivalent with

$$(q-1)^n \varphi_n(x,r) = \sum_{k=0}^{n} (-1)^{n-k} \binom{n}{k} \sum_{j=0}^{k} \begin{bmatrix} k \\ j \end{bmatrix} (q-1)^j q^{r(k-j)} x^j. \qquad (2.15)$$

Comparing coefficients we get (2.12).

## 3. The polynomials $\Phi_n(x,r)$.

Recall that

$$\Phi_n(x,r) = \sum_{k=0}^{n} S(n,k,r) q^{\binom{k}{2}} \left(q^r x\right)^k. \qquad (3.1)$$

**Theorem 3.1**

*The Hankel determinants* $D(n,k) = \det\left(\Phi_{i+j+k}(x)\right)_{i,j=0}^{n-1}$ *are given by*

$$D(n,0) = q^{2\binom{n}{3}+2r\binom{n}{2}} x^{\binom{n}{2}} \prod_{k=0}^{n-1} \left([k]!\left((1-q)x;q\right)_k\right) \qquad (3.2)$$

*and*

$$D(n,1) = D(n,0) \sum_{k=0}^{n} \begin{bmatrix} n \\ k \end{bmatrix} \left(q^{n-1+r} x\right)^k \prod_{j=0}^{n-k-1} [r+j]. \qquad (3.3)$$

Here we set as usual $\prod_{j=0}^{n-1} \left(1-q^j x\right) = (x;q)_n$.

It is easily verified that $\Phi_n(x,r) = q^r x \Phi_{n-1}(qx,r) + x^{-r+1} D x^r \Phi_{n-1}(x,r)$.

This implies

$$\Phi_n(x,r) = x\left(q^r + (q-1)q^r xD + x^{-r}Dx^r\right)\Phi_{n-1}(x,r). \tag{3.4}$$

Let $e(x) = \sum_{n\geq 0} \frac{x^n}{[n]!}$ be the $q-$exponential series (cf. e.g.[1]).

Then it is easily verified that

$x\left(q^r + (q-1)q^r xD + x^{-r}Dx^r\right)x^n = \frac{1}{e(x)x^r}\left(xD\right)x^r e(x)x^n = [r+n]x^n + q^{r+n}x^{n+1}.$

Therefore we get from (3.4)

$$\Phi_n(x,r) = \frac{1}{e(x)x^r}\left(xD\right)^n x^r e(x). \tag{3.5}$$

Since $\frac{1}{x^r}\left(xD\right)^n x^r x^k = [r+k]^n x^k$ formula (3.5) gives immediately

Dobinski's formula

$$\Phi_n(x,r) = \frac{1}{e(x)}\sum_{k\geq 0}\frac{[r+k]^n x^k}{[k]!}. \tag{3.6}$$

Let

$$\left\langle\left\langle x\right\rangle\right\rangle_{r,k} = \prod_{j=0}^{k-1}\frac{x-[r+j]}{q^{r+j}} = q^{-\binom{k}{2}-rk}\left\langle x\right\rangle_{r,k}. \tag{3.7}$$

Define the linear functional $G_r$ by

$$G_r\left(\left\langle\left\langle x\right\rangle\right\rangle_{r,n}\right) = a^n \tag{3.8}$$

and the linear operator $V_r$ by

$$V_r\left(\left\langle\left\langle x\right\rangle\right\rangle_{r,n}\right) = x^n. \tag{3.9}$$

From (1.4) we see that

$$V_r(x^n) = \Phi_n(x,r). \tag{3.10}$$

Then (3.4) implies

$$V_r x V_r^{-1} = x\left(q^r + (q-1)q^r xD + x^{-r}Dx^r\right). \tag{3.11}$$

Consider the polynomials

$$g_n(x,a,r)=\sum_{k=0}^{n}(-a)^k q^{\binom{k}{2}}\begin{bmatrix} n \\ k \end{bmatrix}\langle\langle x\rangle\rangle_{r,n-k}. \tag{3.12}$$

Then $V_r\left(g_n(x,a,r)\right)=p_n(x,a).$

From (3.11) we get

$$\begin{aligned}
&V_r x g_n(x,a,r)=V_r x V_r^{-1}p_n(x,a)=x\left(q^r+(q-1)q^r xD+x^{-r}Dx^r\right)p_n(x,a)\\
&=q^r p_{n+1}(x,a)+q^{n+r}ap_n(x,a)+[r]p_n(x,a)+q^r x[n]p_{n-1}(x,a)+(q^n-1)q^r x^2 p_{n-1}(x,a)\\
&=q^r p_{n+1}(x,a)+q^{n+r}ap_n(x,a)+[r]p_n(x,a)+q^r[n]p_n(x,a)+q^{r+n-1}[n]ap_{n-1}(x,a)\\
&+\left(q^n-1\right)q^r p_{n+1}(x,a)+\left(q^n-1\right)q^r(q^n a+q^{n-1}a)p_n(x,a)+\left(q^n-1\right)q^r q^{2n-2}a^2 p_{n-1}(x,a)\\
&=q^{n+r}p_{n+1}(x,a)+\left([r]+q^r[n]+q^{2n+r}a+q^{2n-1+r}a-q^{n+r-1}a\right)p_n(x,a)\\
&+\left(q^{r+n-1}[n]a+\left(q^n-1\right)q^{2n-2+r}a^2\right)p_{n-1}(x,a)
\end{aligned}$$

Applying $V_r^{-1}$ we get

$$\begin{aligned}
&xg_n(x,a,r)=q^{n+r}g_{n+1}(x,a,r)+\left([r]+q^r[n]+q^{2n+r}a+q^{2n-1+r}a-q^{n-1+r}a\right)g_n(x,a,r)\\
&+\left(q^{r+n-1}[n]a+\left(q^n-1\right)q^{2n-2+r}a^2\right)g_{n-1}(x,a,r)
\end{aligned}$$

In order to get normed polynomials we set

$$H_n(x,a,r)=q^{\binom{n}{2}+rn}g_n(x,a,r). \tag{3.13}$$

Then we have

$$\begin{aligned}
&xH_n(x,a,r)=H_{n+1}(x,a,r)+\left([n+r]+q^{2n+r}a+q^{2n+r-1}a-q^{n+r-1}a\right)H_n(x,a,r)\\
&+q^{2(n-1)+2r}[n]a\left(1+\left(q-1\right)q^{n-1}a\right)H_{n-1}(x,a,r).
\end{aligned} \tag{3.14}$$

It is clear that $G_r\left(H_n(x,a,r)\right)=q^{\binom{n}{2}+rn}p_n(a,a)=0$ for $n>0.$

In the same way as above we get

$$G_r\left(x^n H_n(x,a,r)\right)=\prod_{k=1}^{n}q^{2(k-1)+2r}[k]a\left(1+(q-1)q^{k-1}a\right)=q^{2\binom{n}{2}+2rn}[n]!a^n\prod_{j=0}^{n-1}\left(1+q^j(q-1)a\right).$$

Thus

$$G_r\left(H_n(x,a,r)H_n(x,a,r)\right)=G_r\left(x^nH_n(x,a,r)\right)=q^{2\binom{n}{2}+2rn}[n]!a^n\left((1-q)a;q\right)_n. \qquad (3.15)$$

This immediately implies (3.2). The Hankel determinant (3.3) follows from (3.12) and (1.8).